\documentclass[11pt]{article}

\usepackage{fullwidth}
\usepackage{hyphenat}
\usepackage{color}
\usepackage[color,matrix,arrow]{xy}
\usepackage[top=1.5in, bottom=1.5in, left=1in, right=1in]{geometry}
\usepackage{amsgen}
\usepackage{amsmath}
\usepackage{amstext}
\usepackage{amsbsy}
\usepackage{amsopn}
\usepackage{amsfonts}
\usepackage{amssymb}
\usepackage{eepic}
\usepackage{graphicx}
\usepackage{tikz}
\usepackage{epsf}
\usepackage{pstricks}
\usepackage{mathrsfs}
\usepackage{faktor}
\usepackage{amsthm}
\xyoption{all}

\def\tra#1{\smash{\mathop{\mid\kern
-1pt\joinrel\relbar\joinrel\relbar}\limits^{*}_{#1}}}
\def\longtra#1{\smash{\mathop{\mid\kern
-1pt\joinrel\relbar\joinrel\relbar\joinrel\relbar}\limits^{*}_{#1}}}
\def\vlongtra#1{\smash{\mathop{\mid\kern
-1pt\joinrel\relbar\joinrel\relbar\joinrel\relbar\joinrel\relbar}\limits^{*}_{#1}}}
\def\vvlongtra#1{\smash{\mathop{\mid\kern
-1pt\joinrel\relbar\joinrel\relbar\joinrel\relbar\joinrel\relbar\joinrel\relbar}\limits^{*}_{#1}}}
\def\vvvlongtra#1{\smash{\mathop{\mid\kern
-1pt\joinrel\relbar\joinrel\relbar\joinrel\relbar\joinrel\relbar\joinrel\relbar\joinrel\relbar}\limits^{*}_{#1}}}
\def\etra#1{\smash{\mathop{\mid\kern
-1pt\joinrel\relbar\joinrel\relbar}\limits_{#1}}}

\def\F{{\cal{F}}}
\def\I{{\cal{I}}}

\def\I{{\mathcal{I}}}

\def\bi{\begin{itemize}}
\def\ei{\end{itemize}}
\def\beq{\begin{equation}}
\def\eeq{\end{equation}}

\def\un{\underline}
\def\opn{\operatorname}

\theoremstyle{plain}
\newtheorem{T}{Theorem}[section]
\newcommand{\bt}{\begin{T}}
\newcommand{\et}{\end{T}}
\newcommand{\ftd}{$\square$\end{T}}

\newtheorem{Proposition}[T]{Proposition}
\newcommand{\bp}{\begin{Proposition}}
\newcommand{\ep}{\end{Proposition}}
\newcommand{\fpd}{$\square$\end{Proposition}}

\newtheorem{Lemma}[T]{Lemma}
\newcommand{\bl}{\begin{Lemma}}
\newcommand{\el}{\end{Lemma}}
\newcommand{\fld}{$\square$\end{Lemma}}

\newtheorem{Corol}[T]{Corollary}
\newcommand{\bc}{\begin{Corol}}
\newcommand{\ec}{\end{Corol}}
\newcommand{\fcd}{$\square$\end{Corol}}

\newtheorem{Result}[T]{Result}
\newcommand{\br}{\begin{Result}}
\newcommand{\er}{\end{Result}}
\newcommand{\frd}{$\square$\end{Result}}

\theoremstyle{definition}
\newtheorem{Example}[T]{Example}
\newcommand{\be}{\begin{Example}}
\newcommand{\ee}{\end{Example}}

\newtheorem{Problem}[T]{Problem}
\newcommand{\bq}{\begin{Problem}}
\newcommand{\eq}{\end{Problem}}

\newtheorem{Remark}[T]{Remark}
\newcommand{\brm}{\begin{Remark}}
\newcommand{\erm}{\end{Remark}}

\newtheorem{Definition}[T]{Definition}
\newcommand{\bd}{\begin{Definition}}
\newcommand{\ed}{\end{Definition}}

\newtheorem{Construction}[T]{Construction}
\newcommand{\bco}{\begin{Construction}}
\newcommand{\eco}{\end{Construction}}

\newtheorem{Conjecture}[T]{Conjecture}
\newcommand{\bconj}{\begin{Conjecture}}
\newcommand{\econj}{\end{Conjecture}}

\normalbaselines
\def\abstract#1{\par\bigskip
\begingroup\small
\baselineskip=12truept
\begin{center}ABSTRACT\end{center}
\par\medskip\par\noindent
\null\hfill\hbox{\vbox{\hsize=5truein\noindent#1}}
\hfill\null\par\endgroup\par}

\title{Bases of Permutation Groups and Boolean Representable Simplicial Complexes}
\author{Stuart Margolis and John Rhodes}
\date{\today
\footnote{MCS
Classification: 20B05, 20D05, 20E32}}
\begin{document}
\maketitle

\begin{center}

{\Large Abstract}

\end{center}

\noindent
A base of a permutation group $(X,G)$ is a subset $B$ of $X$ such that the pointwise stabilizer $G_{B}$ is the trivial group. A list $(x_{1},\ldots , x_{k})$ of elements of $X$ is irredundant if each element is not in the pointwise stabilizer of its predecessors. We define a Boolean representable simplicial complex $\mathcal{B}(X,G)$ such that a subset $Y$ of $X$ is independent if and only if some enumeration of its elements is irredundant. In addition $Y$ is a base if and only if its closure is $X$. We give a number of examples and close with a conjecture whose solution leads to a new proof of the Feit-Thompson Theorem.

\section{Introduction}

A base of a permutation group $(X,G)$ is a subset $B$ of $X$ such that the pointwise stabilizer $G_{B}$ of $B$ is trivial. Bases of permutation groups have been studied since the 19th century. They still are a very active and important part of both theoretical and computational group theory.

One reason that bases are important is that the action of $G$ on $X$ is completely determined by its restriction to a base $B$. This allows for the use of bases in algorithmic problems such as computing the order of a permutation group and many other problems.

A base $B$ is {\em irredundant} \cite{CameronBasMat} if there is an enumeration $(x_{1}, \ldots, x_{k})$ of the elements of $B$ such that $x_{i}$ is not in the pointwise stabilizer of $\{x_{1},\ldots, x_{i-1}\}$ for $i=2,\ldots , k$. Bases of permutation groups can be thought of as an analogue of bases of a vector space in the category of permutation groups. However, whereas bases of vector spaces were one of the prime motivators of the theory of matroids, irredundant bases of permutation groups rarely form the bases of a matroid.  In particular elements of an irredundant base need not all have the same size, nor are they invariant under listing them in some order. See Cameron and Fon-der-Flaas \cite{CameronBasMat} for material related to this question. The class of permutation groups that have the property that their bases form a matroid is rare among all permutation groups.

At a conference in 2015 at the University of Durham, Peter Cameron told the first author of this paper that he had been looking for a more general notion of independence that would apply to the case of bases of permutation groups. Margolis responded that the new theory of Boolean Representation \cite{brsc} applied to the theory of bases. Indeed, BRSC give a notion of independence for arbitrary lattices, whereas matroids correspond to the class of geometric lattices. 

Subsequently in December 2015, Cameron gave a lecture entitled ``Beyond Matroids" \cite{BeyondMat} at a Conference at the university of Wellington in honor of Geoff Whittle. In the lecture Cameron wrote that he hadn't yet learned if bases of permutation groups were boolean representable. 

The purpose of this paper is to show that bases and irredundant bases of permutation groups \cite{CameronBasMat} naturally arise in the theory of Boolean representions. More precisely, given a permutation group $(X,G)$ we define a lattice $L=L(X,G)$ with join generators $X$ such that a subset $Y$ of $X$ is dense, that is, its closure is $X$ in the BRSC defined by $L$ \cite{brsc} if and only $Y$ is a base. Furthermore the independent elements in this BRSC are precisely the sets $Y$ such that some enumeration of the elements of $Y$ is a prefix of an irredundant base in the sense of \cite{CameronBasMat}.

We give examples and formulate a conjecture about bases for simple groups. All lattices and simplicial complexes in this paper are assumed to be
finite. Given a set $X$ and $n \geq 0$, we denote by $P_n(X)$
(respectively $P_{\leq n}(X)$) 
 the set of all subsets of $X$ with precisely (respectively at most) 
$n$ elements. All sets in this paper are finite.

\section{Closure Operators, Lattices and Moore Families} \label{colm}

Closure operators on the poset $(P(X), \subseteq)$, lattices and Moore families cryptomorphically define each other. This is well-known material going back to the work of Moore and Ore. See \cite{MooreFam1} for details and an extensive bibliography. We summarize what we need in this paper.

Recall that a closure operator is a function $f:P(X) \rightarrow P(X)$ that is order preserving, idempotent and extensive, the latter meaning that $Y \leq Yf$ for all $Y \subseteq X$. We call a subset $Y=Yf$ in the image of $f$ a closed set.

A subset $\F$ of $P(X)$ is called a {\em Moore family} if $X \in \F$ and if for all $F_{1},F_{2} \in \F$, it follows that $F_{1} \cap F_{2} \in \F$. In other words, a Moore family is a submonoid of the monoid of all subsets of $P(X)$ under intersection. The closed sets $Cl(f)$ of a closure operator $f$ form a Moore Family. Conversely, every Moore family $\F$ defines a closure operator $C(\F)$ defined by letting $(Y)C(\F)=\bigcap_{Z \in \F, Y \subseteq Z} Z$. These operations are inverses of each other and establish a 1-1 correspondence between Moore families and closure operators. 

Furthermore every Moore family $\F$ defines a lattice $L(\F)$ with underlying set $\F$. The meet is given by intersection and the join of two members $F_{1},F_{2}$ of $\F$ is the determined join, that is, the intersection of all members of $\F$ containing $F_{1} \cup F_{2}$. Conversely, if $L$ is a lattice, then the map from $\phi:L \rightarrow (P(L),\cap)$ sending an element $x$ to its principal downset $x^{\downarrow}$ defined by $xf=x^{\downarrow}=\{y \in L| y \leq x\}$ is an isomorphism of meet semilattices to its image which is a Moore family $M(L)$ in $P(L)$. It is easy to check that the lattice defined by $M(L)$ is isomorphic to $L$.

\section{Boolean Representable Simplicial Complexes}

We have established that there is a 1-1 correspondence between closure operators, lattices and Moore families. We use these to define Boolean Representable Simplicial Complexes (BRSC).

A (finite) simplicial complex is a structure of the form $S =
(X,\I)$, where $X$ is a finite nonempty set and $\I \subseteq P(X)$ is nonempty and closed under taking subsets. The elements of $X$ and $\I$ are called respectively {\em points} and {\em independent sets}. 

A maximal independent set is called a {\em basis}. The maximum size of a basis is the {\em rank} of $S$. We say that $S$ is {\em pure} if all its bases have the same size. We say that $S =
(X,\I)$ is {\em simple} if $P_2(X) \subseteq \I$. 

A simplicial complex $M = (X,\I)$ is called a {\em matroid} if it
satisfies the {\em exchange property}:
\bi
\item[(EP)]
For all $I,J \in \I$ with $|I| = |J|+1$, there exists some
  $v \in I\setminus J$ such that $J \cup \{ v \} \in \I$.
\ei

There are many cryptomorphic definitions of matroids \cite{Oxley}. In this paper, since we are concerned with simplicial complexes with various notions of independence, we will always refer to a matroid via its simplicial complex of independent sets as above.

Let $S=(X,\I)$ be a simplicial complex and let $\F$ be a Moore family in $P(X)$ with $\bigcap_{F \in {\mathcal{F}}}F=\square$. We say that $Y\subseteq X$ is a {\em transversal of the
successive differences} for a chain
$$\square = F_{0} \subset F_1 \subset \ldots \subset F_k$$
in $\F$ if $Y$ admits an enumeration $x_1,\ldots , x_k$ such that $x_i \in F_i
\setminus F_{i-1}$ for $i = 1,\ldots,k$. 
We denote by ${\rm Tr}(\F)$ the set of transversals of the
successive differences for chains in $\F$. Alternatively, we can think of a transversal as an ordered list $(x_{1},\ldots , x_{k})$ but we prefer to work with subsets of $X$ and demand that some enumeration of the set $\{x_1 \ldots , x_{k}\}$ is a transversal. We say that a simplicial complex $S = (X,\I)$ is {\em boolean representable}  if $\I\; = {\rm Tr}(\F)$ for some Moore family $\F\, \subseteq P(X)$. 

We  noted in Section \ref{colm} that every Moore Family $\F$ is completely determined by a lattice $L=L(\F)$. $L$ in turn is completely determined by a closure operator $Cl_{L}$. Section 3.6 of \cite{brsc} shows how to formulate BRSC in terms of $L$ and $Cl_{L}$. We use these characterizations in the rest of the paper.


There is a canonical lattice defining a BRSC called its {\em lattice of flats}. If $S=(X,\I)$ is a simplicial complex, then a subset $Y$ of $X$ is called a {\em flat} if 

$$\forall I \in \I \cap P(Y), \forall p \in X\setminus Y, I \cup\{p\} \in \I$$


The set of all flats of $S$ is denoted by 
$L(S)$. Note that $X$ and $\emptyset$ are in $ L(S)$ in all cases, and $L(S)$ is indeed a Moore family that gives $L(S)$ a lattice structure as in Section \ref{colm}.


It follows from \cite[Corollary 5.2.7]{brsc} that a simplicial complex $S = (V,\I)$ is boolean
representable if and only if $\I\, = {\rm Tr}(L(S))$. The associated closure operator on $P(X)$ is defined by
$$Cl(Y) = \cap\{ F \in L(S) \mid Y \subseteq F \}$$
for every $Y\subseteq X$. 

An alternative characterization of BRSC is provided by boolean matrices with a notion of linear independence over the so called super boolean semiring \cite{brsc}. This explains the terminology.

All matroids are boolean representable \cite[Theorem 5.2.10]{brsc}, but the converse is not true. Indeed readers familiar with the theory of matroids will know that the lattice of flats of a matroid $M$ is a geometric lattice and every geometric lattice is the lattice of flats of a matroid \cite{Oxley}. Furthermore independent sets of a matroid are precisely the transversals of its lattice of flats in the sense that we have defined here.
  
However, many lattices that arise naturally in mathematics are not geometric. For example, subgroup lattices of groups or more generally subalgebra lattices of universal algebras are rarely geometric lattices. The same is true for congruence lattices of universal algebras.

The main advantage of BRSC is that they arise from an arbitrary lattice. When we apply this to natural lattices associated to a permutation group we are led precisely to the notion of a base of a permutation group.

\section{Bases of Permutation Groups and Boolean Representable Simplicial Complexes}

Let $G$ be a subgroup of the symmetric group $\operatorname{Sym}(X)$ and $(X,G)$ the corresponding permutation group. If $Y$ is a subset of $X$ we let $G_{Y}$ be the pointwise stabilizer of $Y$. $Y$ is called a {\em base} if $G_{Y}$ is the trivial group. Bases play a crucial role in both the abstract and computational theory of permutation groups. The main theorem of this paper presents bases within the theory of BRSC. For a background on bases, see the 5 lecture series, ``Bases of Permutation Groups" by Professor Tim Burness, University of Bristol \cite{BurnessBases}.

We define a function $f:P(X)\rightarrow P(X)$ by
$Yf =\{z \in X \mid zg =z, \forall g \in G_{Y}\}
=\opn{Stab}(G_{Y})$ where if $H$ is a subgroup of $G$, then $\opn{Stab}(H)=\{x \in X \mid xh=x, \forall h \in H\} $. It is straightforward to show that $f$ is a closure operator. 

Indeed, $f$ comes from the antitone Galois connection $(\phi,\psi)$ between the lattice $(P(X), \subseteq)$ and the subgroup lattice of $G$, $(\opn{Sub}(G), \subseteq)$ defined as follows. If $Y$ is a subset of $X$, let $Y\phi=G_{Y}$ and if $H$ is a subgroup of $G$, let $H\psi=\opn{Stab}(H)$. Clearly $f=\phi\psi$. If $H$ is a subgroup of $G$, then the closure operator $\psi\phi$ applied to $H$ yields the pointwise stabilizer subgroup of $\opn{Stab}(H)$, namely, $G_{\opn{Stab}(H)}$.

The closed sets in the lattice $L=L(X,G)$ 
associated to $f$ are the subsets $Y$ of $X$ such
that for all $x \in X$, $xg=x, \forall g \in
G_{Y}$ implies that $x \in Y$. That is, $Y$ is closed if and only if $Y = \opn{Stab}(G_{Y})$. Recall that a subset $B$
is a base if and only if $G_{B}=\{1\}$. This in turn is true if and only if $\opn{Cl}(B)=X$. Thus we have the following proposition that shows how bases of permutation groups are part of the theory of BRSC.

\bp

Let $(X,G)$ be a permutation group. Then a subset $B$ of $X$ is a base if and only if $\opn{Cl}(B)=X$.

\ep

It is straightforward to see that the closure of the empty set is the set of $x \in X$ such that $xG=\{x\}$. In other words the closure of the empty set is the set of elements whose orbits are trivial. Since these do not effect the questions that we are interested in we exclude trivial orbits in this paper. Therefore in this paper we have $Cl(\emptyset)= \emptyset$. 
We now describe the independent sets in the BRSC $\mathcal{B}(X,G)$ defined by the lattice $L=L(X,G)$ with join generators $X$. Here we send $x \in X$ to $\opn{Stab}(G_{x})$. We have that a set $Y$ is independent if and only if there is an enumeration $(x_{1},x_{2},\ldots , x_{k})$ that forms a transversal of the chain:

$$\emptyset \subset \opn{Stab}(G_{x_{1}}) \subset \opn{Stab}(G_{x_{1}x_{2}})\ldots \subset \opn{Stab}(G_{x_{1}
\ldots x_{k}})$$.

In other words, for each $i=2\ldots k$ $x_{i}$ is not in the pointwise stabilizer of $x_{1}\ldots x_{i-1}$. This means precisely that $Y$ is irrendundant in the sense of \cite{CameronBasMat}, a concept we recall here.

A subset $Y$ of $X$ is {\em irredundant} if there is an enumeration $(x_{1},\ldots , x_{k})$ of elements of $Y$ such that no point $x_{i}$ is stabilized by $G_{\{x_{1},\ldots, x_{i-1}\}}, i=2,\ldots k$. If furthermore, $Y$ is a basis, we call it an {\em irredundant basis}. We have the following theorem.

\bt

Let $(X,G)$ be a permutation group, $L=L(X,G)$ its lattice and $\mathcal{B}(X,G)$ the corresponding Boolean representable simplicial complex. Then a set $Y$ is independent if and only if it is irredundant. In particular, a base is independent if and only if it is an irredundant base.

\et

\proof Let $Y$ be a subset of $X$. Then $Y$ is independent if and only if there is an enumeration $(x_{1},\ldots,x_{k})$ of $Y$ that is a transversal of the chain

$$\emptyset \subset \opn{Stab}(G_{x_{1}}) \subset \opn{Stab}(G_{x_{1}x_{2}})\ldots \subset \opn{Stab}(G_{x_{1}
\ldots x_{k}})$$.

This is true if and only if $(x_{1}, \ldots , x_{k})$ is a transversal of the descending chain of subgroups

$$G \supset G_{x_{1}} \supset G_{\{x_{1},x_{2}\}} \supset \ldots  G_{Y}$$.

One sees this either by definition or by noting the Galois correspondence $(\phi,\psi)$ induces an anti-isomorphism between the closed sets of $(P(X), \subset)$ and $(\opn{Sub}(G), \leq)$ and in particular sends ascending chains in $L(X,G)$ to descending chains in $(\opn{Sub}(G), \leq)$. We therefore see that $Y$ is independent if and only if it is irredundant. $\square$

\brm

Theorem 2.4 of \cite{CameronBasMat} states that if $(X,G)$ is a permutation group then all irredundant bases have the same cardinality if and only if any enumeration of any irredundant base is still an irredundant base if and only if the set of irredundant bases is the set of bases of a matroid.

We point out here that neither of these equivalences are true for the independent sets of a BRSC.  

Let $X=\{1,2,3,4\}$. We write subsets of $X$ by concatenating their elements. Consider the Moore Family $M= \{\emptyset,1,2,3,234,1234\}$. Let $L=L(M)$ be the lattice defined by $M$. Note that $2 \vee 3= 234, 1\vee 2 =1 \vee 3 = 1234$. We take as join generators $\{1,2,3\}$. $(3,2,1)$ is the transversal of the chain $\emptyset \subset 3 \subset 234 \subset 1234$ so that $\{1,2,3\}$ is independent and the BRSC defined by $(L,\{1,2,3\})$, being a simplicial complex, is the uniform matroid on 3 points. However, the enumeration
$(1,2,3)$ is not the transversal of a chain in $L$. Therefore if the BRSC defined by a lattice is a matroid it does not imply that all enumerations of independent sets are transversals.

Example 6.1.4 of \cite{brsc} gives an example of a BRSC $\mathcal{B}$ such that all bases elements have the same length that is, $\mathcal{B}$ is pure, but $\mathcal{B}$ is not a matroid.
\erm




\section{Applications}

Let $(X,G)$ be a permutation group. We will consider the associated permutation group $(P_{2}(X),G)$. We assume that $(X,G)$ is not the right regular representation of the cyclic group of order 2 so to ensure that $(P_{2}(X),G)$ is a permutation group. The action is $\{x,y\}g=\{xg,yg\}$. For ease of notation we write $ab$ for $\{a,b\}$ so that $ab=ba, a\neq b \in X$. 

If $H$ is a subset of $G$ and $Z$ is a subset of $X$, then $Z$ is $H$-invariant if $ZH=Z$. This does not necessarily imply that $H$ fixes $Z$ pointwise. The $H$-invariant subsets of $X$ form a
Boolean algebra under union, intersection and complement of subsets. We give two classes of examples.

\be

Let $\opn{Sym}(n)$ be the symmetric group on the set $\un{n}=\{1, \ldots, n\}$.
Consider $(P_{2}(\un{n}),\opn{Sym}(n))$ for $n>2$.

Assume that $n$ is divisible by 3. Then $\{12,23,45,56,\ldots,(n-2)(n-1), (n-1)(n)\}$ is a basis. Indeed if $g \in \opn{Sym}(n)$ fixes $12$ and $23$, then $2g=2$ and thus $g$ fixes $\{1,2,3\}$ pointwise. Now continue by induction.

A small modification in case $n$ is congruent to 1 or 2 modulo 3 shows that $(P_{2}(\un{n}),\opn{Sym}(n))$ has a basis of size about $(2/3)n$.

\ee

\be \label{odd}

Let $G$ be a group of odd order and consider $(P_{2}(X),G)$ where $(X,G)$ is a permutation group with $|X|=n$. Then $12,34,56\ldots$ is a basis of size $n/2$ if $n$ is even and $(n-1)/2$ if $n$ is odd.

Indeed if $g \in G$ acts invariantly on $12$ it fixes $12$ pointwise since $G$ has odd order. Now continue by induction.

\ee

The last example leads to the following conjecture for simple non-abelian groups.

\bconj

Let $G$ be a simple non-abelian group. Let $n'=n/2$ if $n$ is even and $n'=(n-1)/2$ if $n$ is odd. Then there exists a permutation group $(X,G)$ such that every base $B$ for
$(P_{2}(X),G)$ has $|B|$ greater than $n'$ where $n=|X|$.

\econj

We note that if this conjecture is true, then it follows from Example \ref{odd} that the order of $G$ is even and the Feit-Thompson Theorem follows from the conjecture.





\bibliography{stubib}
\bibliographystyle{abbrv}

{(Stuart Margolis) Department of Mathematics, Bar Ilan University, Ramat Gan 52900, Israel

{\it Email address}: \; \texttt{margolis@math.biu.ac.il}}

\medskip

{(John Rhodes) Department of Mathematics, University of California, Berkeley, CA 94720, U.S.A.

{\it Email address}: \; \texttt{rhodes@math.berkeley.edu, blvdbastille@gmail.com}}

\end{document}